\def\patterns#1{}
\input plain
\input amstex
\documentstyle{amsppt}  
\pagewidth{12.5cm}\pageheight{19cm}\magnification\magstep1
\topmatter
\title Canonical bases in Lie theory and total positivity\endtitle
\author G. Lusztig\endauthor
\abstract{This represents a talk given at the International
Conference for Basic Science, July 2025. We review the theory of
canonical bases of quantum groups and its
relation with the theory of total positivity.}\endabstract

\address{Department of Mathematics, M.I.T., Cambridge, MA 02139}\endaddress
\endtopmatter   
\document

\define\mpb{\medpagebreak}

\define\si{\sim}

\define\sqc{\sqcup}

\define\op{\oplus}
   
\define\part{\partial}

\define\imp{\implies}

\define\m{\mapsto}
\define\do{\dots}

\define\lra{\leftrightarrow}

\define\nl{\newline}
\redefine\i{^{-1}}

\define\Hom{\text{\rm Hom}}
\define\End{\text{\rm End}}

\define\tr{\text{\rm tr}}

\define\card{\text{\rm card}}

\define\a{\alpha}

\redefine\d{\delta}
\define\e{\epsilon}

\redefine\l{\lambda}

\redefine\L{\Lambda}

\define\boc{\bold c}

\define\ii{\bold i}

\define\kk{\bold k}

\define\BB{\bold B}
\define\CC{\bold C}

\define\LL{\bold L}

\define\NN{\bold N}

\define\QQ{\bold Q}
\define\RR{\bold R}

\define\UU{\bold U}
\define\VV{\bold V}

\define\ZZ{\bold Z}

\define\cb{\Cal B}

\define\ff{\frak f}

\define\fU{\frak U}

\define\fZ{\frak Z}

\define\bul{\bullet}

\head Introduction\endhead
\subhead 0.1\endsubhead
The 1888-1890 paper \cite{K} by W. Killing on the classification of (finite
dimensional) simple Lie algebras over $\CC$ has been called
``the greatest mathematical paper of all time''
(A. J. Coleman \cite{Co}).

The arguments of \cite{K} were made more rigorous by
\'E. Cartan (in 1894, see \cite{Ca}),
who also classified the (finite dimensional) irreducible
representations.

The simple Lie algebras are of the following types:

$A_n,D_n,E_6,E_7,E_8$  (simply laced)

$B_n,C_n,F_4,G_2$ (non-simply laced,related to simply laced ones by folding).

We will focus on simply laced types.

\subhead 0.2\endsubhead
The results of Killing and Cartan can be interpreted as saying that
the simply laced simple Lie algebras are in bijection with a certain
set of graphs, called Coxeter graphs.
(These graphs have a simple characterization, due to Coxeter, see
\cite{L23b}.)

\subhead 0.3\endsubhead
Here are some developments after the work of Killing-Cartan.

(1) The irreducible representations defined by
Cartan have character/dimension described by H. Weyl in 1925.

(2) The simple Lie groups over $\CC$ have analogues where $\CC$
is replaced by any algebraically closed field (C. Chevalley, 1955).

(3) In 1985, Drinfeld, Jimbo defined quantum groups, some
algebras over $\QQ(v)$ ($v$ is an indeterminate) which are
deformations of the enveloping algebra of a simple Lie algebra.

(4) (A mixture of (3),(2).) In \cite{L89} I defined a form over
$\QQ[v,v^{-1}]$ of the quantum group in (3) and used it to define a
version of the quantum group in which the parameter $v$ is a root
of $1$. This gave new information on the representation theory of
the groups in (2) in characteristic $>1$. 

(5) In \cite{L90} I observed that the
theory of quantum groups leads to some new,
extremely rigid structures in which the 
objects of the theory are provided with  canonical bases
with rather remarkable properties; in particular,
the structure constants with respect to the canonical
basis are in $\NN[v,v^{-1}]$.
These specialize for $v=1$ to canonical bases
for the objects in the classical theory (such
as the irreducible representations constructed
by Cartan) in which the structure constants
are in $\NN$. See \cite{L25} for a history
of canonical bases.

(6) In \cite{L94} I observed that simple
Lie groups over $\CC$ have semigroup
analogues where $\CC$ is replaced by certain
semifields $K$ (for example $K=\RR_{>0}$).
The case of type $A$ with $K=\RR_{>0}$ was first
considered by Schoenberg and Gantmacher-Krein around 1930.

In this paper we will focus on (5),(6). Our main theme is
that (5),(6) are closely related. We will show several instances
where something about (5) implies something about (6) (we then write
$(5)\imp(6)$) and several instances where something about (6) implies
something about (5) (we then write $(6)\imp(5)$).

\subhead 0.4. Notation\endsubhead
In this paper we fix a Coxeter graph with set of vertices $I$.
For $i,j$ in $I$ let $n_{ij}$ be the number of edges
that join $i,j$. We have $n_{ij}\le 1$ for all $i,j$ and
$n_{ii}=0$ for all $i$.

\mpb

Let $V$ be the $\QQ$-vector space with
basis $\{\a_i;i\in I\}$. For $i\in I$ define $s_i:V@>>>V$
      (linear) by

$s_i(\a_i)=-\a_i$

$s_i(\a_j)=\a_i+\a_j$ if $n_{ij}=1$

$s_i(\a_j)=\a_j$ if $i\ne j,n_{ij}=0$.

Let $W$ be the subgroup of $GL(V)$ generated
by $\{s_i;i\in I\}$; a finite (Weyl) group. Let

$R=\{x\in V;x=w(\a_i) \text{ for some }i\in I,w\in W\}$ (roots).

We have $R=R^+\sqcup(-R^+)$ where
$R^+=R\cap\sum_{i\in I}\NN\a_i$.

There is a unique $w_0\in W$ such that $w_0(R^+)=-R^+$.

\head Contents\endhead
1. The simple Lie algebra and the simple algebraic group
attached to $I$.

2. The semigroup $\fU(K)$ for a semifield $K$. 

3. The involution $\phi_K:\fU(K)@>>>\fU(K)$ for a semifield $K$.

4. The canonical basis $\BB$.

5. The canonical basis $\BB_d$.

6. Using quivers.

7. The semigroup $G(\RR_{>0})$.

8. The semigroup $G(K)$ for a semifield $K$.

\head 1. The simple Lie algebra and the simple algebraic group
attached to $I$\endhead
\subhead 1.1\endsubhead
Let $\kk$ be an algebraically closed field.
Let $L_\kk$ be the $\kk$-vector space with basis

$\{X_\a,\a\in R\}\sqcup\{t_i,i\in I\}$.

For $i\in I,\e=\pm1$ define a linear
map ${}^\e\underline i:L_\kk@>>>L_\kk$ by

$X_\a\m0$ if $\a\in R,\a+\e\a_i\notin R\cup 0$,

$X_\a\m X_{\a+\e\a_i}$ if
                 $\a\in R,\a+\e\a_i\in R$,

$X_\a\m t_i$ if $\a\in R,\a+\e\a_i=0$,

$t_i\m 2X_{\e\a_i}$,

$t_j\m n_{ij}X_{\e\a_i}$ if $i\ne j$.

When $\kk=\CC$, the Lie subalgebra of $End(V)$
generated by ${}^\e\underline i,i\in I,\e=\pm1$
is the simple Lie algebra 
      attached by Killing to our Coxeter graph.

      (This is not the standard definition
      using generators and Serre relations; we do not use relations.
      We use the theory of  canonical bases or rather a precursor
of that theory in \cite{L90a}.)           

\subhead 1.2\endsubhead
      For general $\kk$ and for $i\in I,\e=\pm1,\l\in\kk$
     define an (invertible ) linear map
      ${}^\e i^\l:L_\kk@>>>L_\kk$ by

      $X_\a\m X_\a$ if
      $\a\in R,\a+\e\a_i\notin R\cup 0$,

$X_\a\m X_\a+\l X_{\a+\e\a_i}$ if $\a\in R,\a+\e\a_i\in R$,

$X_\a\m X_\a+\l t_i +\l^2X_{\e\a_i}$ if $\a\in R,\a+\e\a_i=0$,

$t_i\m t_i+2\l X_{\e\a_i}$,

$t_j\m t_j+\l X_{\e\a_i}$ if $n_{ij}=1$,

$t_j\m t_j$ if $i\ne j$, $n_{ij}=0$.

The subgroup of $GL(V)$ generated by
      ${}^\e i^\l,i\in I,\e=\pm1,\l\in\kk$
is the simple algebraic group $G(\kk)$ 
attached by Chevalley (1955) to $\kk$ and our Coxeter graph.

(This is not the standard definition.
It uses the theory of canonical bases
or rather a precursor of that theory in \cite{L90a}.)           

When $\kk=\CC$ we have
      ${}^\e i^\l=\exp(\l{}^\e\underline i)$.

\head 2. The semigroup $\fU(K)$ for a semifield $K$\endhead
\subhead 2.1\endsubhead
Let $K$ be a semifield, that is a set with
two commutative, associative operations $+,\times$ such that
      
$(a+b)c=ac+bc$ for $a,b,c$ in $K$ and
$K$ is a group with respect to $\times$.

Here are some examples.

$K=\RR_{>0}$ with the usual $+,\times$.

$K=\ZZ$ with $a+b=\min(a,b)$, $ab=\text{ usual }a+b$.

$K=\{1\}$ with $1+1=1,11=1$.

\subhead 2.2\endsubhead
Let $\fU^\bul(K)$ be the semigroup with generators
$i^a,i\in I,a\in K$ and relations
$$i^aj^bi^c=j^{a'}i^{b'}j^{c'}$$
if $n_{ij}=1$ and $a,b,c,a',b',c'$ in $K$
satisfy $a'b'=bc,b'c'=ab$, $b'=a+c$ or equivalently $b=a'+c'$,
$$i^aj^b=j^bi^a$$
if $i\ne j,n_{ij}=0$ and $a,b$ in $K$,
$$i^ai^b=i^{a+b}$$
if $a,b$ in $K$.

When $K=\RR_{>0}$ this is the same as the
subsemigroup of $G(\CC)$ generated by 
${}^1i^\l,i\in I,\l\in\RR_{\ge0}$. (This is an instance of (6) in 0.3.)

Let $*:\fU^\bul(K)@>>>\fU^\bul(K)$ be the semigroup antiautomorphism
such that $*(i^a)=i^a$ for all $i\in I,a\in K$.

\subhead 2.3\endsubhead
Let $\fZ$ be the set of all sequences
$\ii=(i_1,i_2,\dots,i_\nu)$ in $I$ such that
$\nu=\sharp(R^+)$,
$w_0=s_{i_1}s_{i_2}\dots s_{i_\nu}$.

For any $\ii=(i_1,i_2,\dots,i_\nu)\in\fZ$, the map
$x_{\ii}:K^\nu@>>>\fU^\bul(K)$,
$$(c_1,c_2,\dots,c_\nu)\m i_1^{c_1}i_2^{c_2}\dots i_\nu^{c_\nu}$$
is injective and its image is independent
of $\ii$; this image is denoted by $\fU(K)$.
It is a subsemigroup of $\fU^\bul(K)$.

Note that $*:\fU^\bul(K)@>>>\fU^\bul(K)$ restricts to a semigroup
antiautomorphism of $\fU(K)$.

When $K$ is the semifield $\ZZ$ and $\ii\in\fZ$,
then the subset $x_{\ii}(\NN^\nu)$ of $\fU(\ZZ)$ is
independent of $\ii$; we denote it by $\fU(\NN)$.
This is a subsemigroup of $\fU(\ZZ)$.

\head 3. The involution $\phi_K:\fU(K)@>>>\fU(K)$ for
a semifield $K$\endhead
\subhead 3.1\endsubhead
Let $K$ be a semifield.
For $i\in I,a\in K$ there is a well defined
bijection $\tau_{i,a}:\fU(K)@>>>\fU(K)$ such that
$$\tau_{i,a}(i_1^{c_1}i_2^{c_2}\dots i_\nu^{c_\nu})=
i_1^{ac_1}i_2^{c_2}\dots i_\nu^{c_\nu}$$
whenever
$\ii=(i_1,i_2,\dots,i_\nu)\in\fZ$ satisfies $i_1=i$ and
$\boc=(c_1,c_2,\dots,c_\nu)\in K^\nu$.

Let $\ii\in\fZ$. For $k\in[1,\nu]$ we write (in $V$):
$$s_{i_1}s_{i_2}\dots s_{i_{k-1}}(\a_{i_k})=
\sum_{i\in I}r_{i,k}\a_i$$
where $r_{i,k}\in\NN$. Define $q:I@>>>K$ by
$$q(i)=\sum_{k\in[1,\nu]}r_{i,k}\in\ZZ_{>0}\subset K.$$

For $p:I@>>>K$ we define
$$\boc_p=(c_1,c_2,\dots,c_\nu)=
(\sum_ir_{i,1}p_{i_1},\sum_ir_{i,2}p_{i_2},\dots,
\sum_ir_{i,\nu}p_{i_\nu})\in K^\nu.$$
Then
$$u(p)=i_1^{c_1}i_2^{c_2}\dots i_\nu^{c_\nu}\in\fU(K)$$
is independent of the choice of $\ii$, see \cite{L94}.

For $i\in I$ we define $i^!\in I$ by $w_0(\a_i)=-\a_{i^!}$.

We now state the following result.
\proclaim{Theorem 3.2}
There is a unique bijection $\phi_K:\fU(K)@>>>\fU(K)$ such that

(a) $\tau_{i,a}\phi_K=\phi_K\tau_{i^!,a^{-1}}$
for any $i\in I,a\in K$ and

(b) $\phi_K(1)=u(q^{-1})$ where $1:I@>>>K$ is constant with value $1$.
\endproclaim
When $K=\RR_{>0}$ we define $\phi_K$ to be $\phi$ in 7.4.
Then (a),(b) are satisfied by \cite{L23b, \S11, A1, A2}. We then define
$\phi_K$ for any semifield $K$ by the arguments in \cite{L23a, 4.7};
then (a),(b)
continue to hold.
This proves the existence statement. We now prove uniqueness.
We regard $\fU_K$ as the set of vertices of a graph as in
\cite{L23b, \S9}. To prove uniqueness it is enough to show
that this graph is connected. It is known that we can find a
surjective homomorphism of semifields $K'@>>>K$ where $K'$ is
contained in the set of nonzero elements of a field. For such $K'$
the connectedness of the graph above is proved in \cite{L24}.
This completes the proof.

Note that we have $\phi_K^2=1$.

\subhead 3.3\endsubhead
For example,

if $I=\{i\}$ then $\phi_K(i^a)=i^{1/a}$ for $a\in K$;

if $I=\{i,j\}$, $n_{ij}=1$ then
$$\phi_K(i^aj^bi^c)=i^{a/c(a+c)}j^{(a+c)/ab}i^{1/(a+c)}
=j^{c/ab}i^{1/c}j^{1/b}$$
for $a,b,c$ in $K$.

\subhead 3.4\endsubhead
For any $p:I@>>>K$, let
$S_p:\fU(K)@>>>\fU(K)$ be the unique bijection such that
$$i_1^{c_1}i_2^{c_2}\dots i_\nu^{c_\nu}\m
 i_1^{c_1p(i_1)}i_2^{c_2p(i_2)}\dots i_\nu^{c_\nu p(i_\nu)}$$
for any $\ii\in\fZ,\boc\in K^\nu$.

\head 4. The canonical basis $\BB$ \endhead
\subhead 4.1\endsubhead
Let $\ff$ be the $\QQ(v)$-algebra with generators
$E_i,i\in I$ and relations

$E_i^2E_j-(v+v^{-1})E_iE_jE_i+E_jE_i^2=0$
if $n_{ij}=1$,

$E_iE_j=E_jE_i$ if $i\ne j,n_{ij}=0$.

(This is the $+$ part of the Drinfeld-Jimbo
quantum group. There is also an alternative definition in which
no relations are used, see \cite{L93}.)

There is a unique $\QQ$-algebra $\QQ$-algebra
isomorphism $\bar{}:\ff@>>>\ff$ such that 

$\bar{E_i}=E_i$, $\bar{v^n}=v^{-n}$ for $i\in I,n\in\ZZ$.

For $x\in\ff$ and $n\in\NN$ we set
$x^{(n)}=x^n/[n]!$
where
$$[n]!=\prod_{h=1}^n(v^h-v^{-h})/(v-v^{-1}).$$

\subhead 4.2\endsubhead
For $w\in W,i\in I$ such that $w(\a_i)\in R^+$
we define $T'_{w,-1}(E_i)\in\ff$ in terms of a braid
group action on the full quantum group as in \cite{L93, 40.1.2};
this is a ``v-analogue of root vectors''. 

For example, if $w=s_j$ and $s_j(\a_i)=\a_i+\a_j$ then
$T'_{w,-1}(E_i)=E_iE_j-v^{-1}E_jE_i$.

\subhead 4.3\endsubhead
Let $\ii=(i_1,i_2,\dots,i_\nu)\in\fZ$.
For any $\boc=(c_1,c_2,\dots,c_\nu)\in\NN^\nu$ we set
$$E_{\ii,\boc}
=E_{i_1}^{(c_1)}T'_{s_{i_1},-1}(E_{i_2})^{(c_2)}\dots
T'_{s_{i_1}s_{i_2}\dots s_{i_{\nu-1}},-1}(E_{i_\nu})^{(c_\nu)}
\in\ff.$$

One can show that for any $\ii\in\fZ$,
$$B_\ii=\{E_{\ii,\boc};\boc\in\NN^\nu\}$$
is a $\QQ(v)$-basis of $\ff$.

 \proclaim{Theorem 4.4, see \cite{L90}}

(i) The $\ZZ[v^{-1}]$-submodule $\Cal L$
spanned by $B_\ii$ is independent of $\ii$.

(ii) The image of $B_\ii$ under the obvious homomorphism
$\pi:\Cal L@>>>\Cal L/v^{-1}\Cal L$ is a $\ZZ$-basis $\beta$ of
$\Cal L/v^{-1}\Cal L$ independent of $\ii$. 

(iii) $\pi$ restricts to an isomorphism of $\ZZ$-modules
$\Cal L\cap\bar\Cal L@>>>\Cal L/v^{-1}\Cal L$.

(iv) For $b\in\beta$ define $\tilde b\in\Cal L\cap\bar\Cal L$ by
$\pi(\tilde b)=b$.
We have $\bar{\tilde b}=\tilde b$.

(v) $\BB=\{\tilde b;b\in\beta\}$ is a $\QQ(v)$-basis of $\ff$,
a $\ZZ[v^{-1}]$-basis of $\Cal L$
and a $\ZZ$-basis of $\Cal L\cap\bar\Cal L$.

(vi) We have a bijection $\fU(\NN)@>\si>>\BB$,
$$i_1^{c_1}i_2^{c_2}\dots i_\nu^{c_\nu}\m\widetilde{\pi(E_{\ii,\boc})}.
$$
\endproclaim
Note that (vi) is an instance of $(6)\imp(5)$ and of $(5)\imp(6)$, see
0.3.

\subhead 4.5\endsubhead
If $I=\{i\}$ (type $A_1$) then $\BB$ consists of
$E_i^{(a)}$ for various $a\in\NN$.

If $I=\{i,j\}$ with $n_{ij}=1$ (type $A_2$)
then $\BB$ consists of
$$E_i^{(a)}E_j^{(b)}E_i^{(c)},E_j^{(a)}E_i^{(b)}E_j^{(c)}$$
for various $a,b,c$ in $\NN$ such that $b\ge a+c$ where the elements
$$E_i^{(a)}E_j^{(a+c)}E_i^{(c)}=E_j^{(c)}E_i^{(a+c)}E_j^{(a)}$$
are considered only once.

Assume that $I=\{i,j,k\}$ with $n_{ij}=1,n_{jk}=1, n_{ik}=0$.
(Type $A_3$). Then 
$$E_j^{(2)}E_iE_kE_j-E_j^{(3)}E_iE_k=E_jE_iE_kE_j^{(2)}-
E_iE_kE_j^{(3)}\in\BB.$$
In \cite{X} all elements of $\BB$ in type $A_3$ are described.

In \cite{A} an algorithm to compute $\BB$ is given.

\head 5. The canonical basis $\BB_d$\endhead
\subhead 5.1\endsubhead
Let $d:I@>>>\NN$ and let
$$J_d=\sum_{i\in I}\ff E_i^{d(i)+1}\subset\ff.$$
Let $\LL_d=\ff/J_d$.
It is known that $\LL_d$ is an irreducible
representation of the quantum group $\UU$
with generators $E_i,F_i,K_i$ attached to $I$ in which

$E_i$ acts by left multiplication,

$F_i$ maps $1$ to $0$,

and $K_i$ maps $1$ to $v^{-d(i)}$.

\proclaim{Theorem 5.2, see\cite{L90}}
$\BB\cap J_d$ is a $\QQ(v)$-basis
of $J_d$ hence the obvious map $\ff@>>>\LL_d$
restricts to a bijection of $\BB_d:=\BB-(\BB\cap J_d)$
onto a $\QQ(v)$-basis of $\LL_d$ (denoted again by $\BB_d$).
\endproclaim

\subhead 5.3\endsubhead
Specializing $v$ to $1$ we see that any
irreducible representation (defined by Cartan) of the Lie
algebra attached to $I$ has a canonical basis.

In this basis the action of ${}^1\underline i,
{}^{-1}\underline  i$ is by matrices with entries in $\NN$.

(The result for ${}^{-1}\underline i$ is reduced to the result for
${}^1\underline i$ by  5.5(a).)

\subhead 5.4\endsubhead
According to \cite{L23a}:

The bijection $\fU(\NN)@>>>\BB$ restricts to a bijection
of the subset
$$\fU^d(\NN)=\{A\in\fU(\NN);S_d\phi_\ZZ(A)\in\fU(\NN)\}$$
of $\fU(\NN)$ with the subset $\BB_d$ of $\BB$.
This is an instance of $(6)\imp(5)$ and of $(5)\imp(6)$, see 0.3.

\subhead 5.5\endsubhead
Let $\tilde \LL_d$ be the vector space $\LL_d$ with a new
$\UU$-module structure such that

$E_i$ acts on $\tilde\LL_d$ as $F_i$ acts on $\LL_d$,

$F_i$ acts on $\tilde\LL_d$ as $E_i$ acts on $\LL_d$,

$K_i$ acts on $\tilde\LL_d$ as $K_i^{-1}$ acts on $\LL_d$.

Define $d':I@>>>\NN$ by $d'(i)=d(j)$ where $\a_j=-w_0(\a_i)$.

The following result appears in \cite{L90b}.

(a) There is a unique isomorphism of $\UU$-modules
$\LL_{d'}@>>>\tilde\LL_d$ which carries $\BB_{d'}$ onto $\BB_d$. 

Now $E_i\m E_{i^!}$ defines an algebra automorphism of $\ff$
which restricts to a bijection $\BB@>>>\BB$ and to a bijection
$\BB_d@>>>\BB_{d'}$. Composing this with the bijection
$\BB_{d'}@>>>\BB_d$ above we obtain a bijection $\BB_d@>>>\BB_d$.

\subhead 5.6\endsubhead
Now $A\m S_d\phi_\ZZ(A)$ is an involution of the set $\fU(\ZZ)$.
It restricts to an involution of $\fU^d(\NN)$. This can
be viewed as an involution of $\BB_d$.

This coincides with the involution of $\BB_d$ described in 5.5.

\subhead 5.7\endsubhead
Note that Weyl's dimension formula expresses $\dim\LL_d$ as a quotient
of two integers (the result being not obviously an integer).

Our equality $\dim\LL_d=\sharp(\fU^d(\NN))$ expresses directly
$\dim\LL_d$ as a cardinal of an explicit finite set.

\head 6. Using quivers\endhead
\subhead 6.1\endsubhead
Let $H$ be the set of ordered pairs $h=(i,j)$ of elements in $I$
such that $n_{ij}=1$. For such $h$ we set $h'=i,h''=j,\bar h=(j,i)$.
We fix a function $e:H@>>>\CC^*$ such that $e(h)+e(\bar h)=0$ for any
$h\in H$.
Let $\Cal C$ be the category of finite dimmensional $I$-graded
$\CC$-vector spaces $\VV=\op_{i\in I}\VV_i$.
For $\VV\in\Cal C$ we set
$$E_\VV=\oplus_{h\in H}\Hom(\VV_{H'},\VV_{h''}).$$
We define a nondegenerate symplectic form on $E_\VV$ by
$$(x,y)=\sum_{h\in H}e(h)\tr(x_hy_{\bar h}).$$
Following Gelfand and Ponomarev, for any $i\in I$ we define
$\Psi_i:E_\VV@>>>\End(\VV_i)$ by
$$\Psi_i(x)=\sum_{h\in H;h''=i}e(h)x_hx_{\bar h}.$$
We set
$$\L_\VV=\{x\in E_\VV;\Psi_i(x)=0 \text{ for all }i\in I\}.$$
By  \cite{L91, 12.3(a) ,14.1(a)},
$\L_\VV$ is a closed subvariety of $E_\VV$ of
pure dimension $\dim E_\VV/2$. Moreover it is Lagrangian.
Let $Irr\L_\VV$ be the set of irreducible components of $\L_\VV$. 

Restricting the bijection obtained by combining
\cite{L92, 4.16(a), 6.16(a)} (in the affine case)
we obtain a bijection
$$\BB@>>>\sqc_\VV Irr(\L_\VV).$$
(Here $\VV$ runs over the set of isomorphism classes of objects
in $\Cal C$.)  Composing this with the bijection $\fU(\NN)@>>>\BB$
we obtain a bijection

$\fU(\NN)@>>>\sqc_\VV Irr(\L_\VV)$.
\nl
Let $\xi_X\in\fU(\NN)$ be the element corresponding to 
$X\in Irr(\L_\VV)$ under this bijection.

\subhead 6.2\endsubhead
We now describe explicitly the element $\xi_X$.
There is a unique partition $I=I_0\sqc I_1$ such that $n_{ij}=0$ for any $i\ne j$ in $I_0$ and also for any $i\ne j$ in $I_1$.
For $\d\in\{0,1\}$ let $H_\d$ be the set of all
$(i,j)\in H$ such that $j\in I_\d$ (hence $i\in I_{1-\d}$).
We have $E_\VV=E_\VV^0\op E_\VV^1$ where
$$E_\VV^\d=\{x\in E_\VV;x_h=0 \text{ if }h\in H_{1-\d}\}.$$

Let $r_0=\card(I_0),r_2=\card(I_1)$.
For $\d\in\{0,1\}$ we consider a sequence $\ii_\d=(i_{1,\d},i_{2,\d},
\do,i_{\nu,\d})$
in $I$ in which
the first $r_\d$ terms are the elements of $I_\d$ in some order,
the next $r_{1-\d}$ terms are the elements of $I_{1-\d}$ in some order,
the next $r_\d$ terms are the elements of $I_\d$ in some order,
the next $r_{1-\d}$ terms are the elements of $I_{1-\d}$ in some order,
etc.

(We have $\nu=r_\d+r_{1-\d}+r_\d+r_{1-\d}+\dots$, a sum of
$2\nu/\card(I)$ terms.) It is known that $\ii_\d\in\fZ$.
Note that $\ii_\d$ is adapted to the orientation $H_\d$ in the
sense of \cite{L90, 4.7}.

Let $G_\VV=\prod_{i\in I}\VV_i$. This group acts on $E_\VV$
by $(g_i):(x_h)\m(y_h)$ where $y_h=g_{h''}x_hg_{h'}\i$.
This restricts to an action of $G_\VV$ on $E^\d_\VV$ for
$\d\in\{0,1\}$.
Let $\Cal F_\VV$ be the set of all functions $f:R^+@>>>\NN$
such that 
$$\sum_{\a\in R^+}f(\a)=\sum_{i\in I}\dim(\VV_i)\a_i.$$

Using a theorem of Gabriel \cite{G},
proved earlier by Brauer (unpublished)  for type $A,D$,
we see that the set of $G_\VV$-orbits on $E^\d_\VV$
is in canonical bijection with the set $\Cal F_\VV$ consisting
of  all functions $f:R^+@>>>\NN$ such that 
$$\sum_{\a\in R^+}f(\a)=\sum_{i\in I}\dim(\VV_i)\a_i.$$
(This bijection has been used by Ringel \cite{R} to relate
quantum groups with Hall algebras.)
We denote by $\Cal O^\d_f$ the $G_\VV$-orbit corresponding to $f$.

Let $\d\in\{0,1\}$. As in \cite{L91,12.8} we can identify
$E_\VV$ with the cotangent bundle of $E^\d_\VV$. If
$X\in Irr(\L_\VV)$ then by \cite{L91,14.3}, $X$ is the
closure of the conormal bundle of a well defined $G_\VV$-orbit
$\Cal O^\d_X$ in $E^\d_\VV$. Let $f^\d_X:R^+@>>>\NN$ be the
function in $\Cal F_\VV$ corresponding to this orbit as above.
We set
$$\xi_X^\d=
(i_{1,\d}^{c_1}i_{2,\d}^{c_2}\do i_{\nu,\d}^{c_\nu})\in\fU(\NN)$$
where

$c_1=f^\d_X(\a_{i_{1,\d}})$,

$c_2=f^\d_X(s_{i_{1,\d}}(\a_{i_{2,\d}}))$,

$c_3=f^\d_X(s_{i_{1,\d}}s_{i_{2,\d}}(\a_{i_{3,\d}}))$,

etc.

One can show that $\xi_X=\xi_X^\d$.

\subhead 6.3\endsubhead
The consideration of the intersection cohomology sheaves of the closures
of $G_\VV$-orbits on $E^\d_\VV$ (see 6.2) leads to
an alternative (geometric) definition of $\BB$ (see \cite{L90}). 

This provides very strong positivity properties of $\BB$ which are
not seen from the purely algebraic approach.

\head 7. The semigroup $G_{>0}$\endhead
\subhead 7.1\endsubhead
There is a unique automorphism $\Omega:G(\CC)@>>>G(\CC)$ such that
$\Omega({}^\e i^\l)={}^{-\e}i^\l)$ for $i\in I,\e=\pm1,\l\in\CC$.
We have $\Omega^2=1$.

Let $U^+$ (resp. $U^-$) be the subgroup of $G(\CC)$
generated by the elements  ${}^1i^\l$ (resp. ${}^{-1}i^\l$) with
$i\in I,\l\in\CC$.
Let $B^+$ (resp. $B^-$) be the normalizer of $U^+$ (resp. $U^-$)
in $G(\CC)$ (a Borel subgroup). 
Let $U^+_{\ge0}$ (resp. $U^-_{\ge0})$ be the subsemigroup of
$U^+$ (resp. $U^-$) generated by
the elements ${}^1i^\l$ (resp, ${}^{-1}i^\l$) with
$i\in I,\l\in\RR_{\ge}$.
We write $\fU_{>0}$ instead of $\fU(\RR_{>0})$.
We have a semigroup isomorphism
$\fU^\bul(\RR_{>0})@>>>U^+_{\ge0}$, $u\m u^+$,
(resp. $\fU^\bul(\RR_{>0})@>>>U^-_{\ge0}$, $u\m u^-$)
such that $i^\l\m {}^1i^\l$ (resp. $i^\l\m {}^{-1}i^\l$)
for $i\in I,\l\in\RR_{\ge0}$.
Let $U^+_{>0}=\{u^+;u\in\fU_{>0}\}$,
$U^-_{>0}=\{u^-;u\in\fU_{>0}\}$.
Let $T=B^+\cap B^-$ (a maximal torus).
Let $T_{>0}$ be the subgroup of $T$ generated by the
elements $\chi(\l)$ for various homomorphisms of algebraic groups
$\chi:\CC^*@>>>T$ and various $\l\in\RR_{>0}$.
We have $U^+_{>0}T_{>0}U^-_{>0}=U^-_{>0}T_{>0}U^+_{>0}$.
This is a subsemigroup of $G(\CC)$ denoted by $G_{>0}$.

\subhead 7.2\endsubhead
Let $\cb$ be the variety of Borel subgroups of $G(\CC)$. The maps

$c^+:\fU_{>0}@>>>\cb, u\m u^+B^-(u^+)\i$,
$c^-:\fU_{>0}@>>>\cb, u\m u^-B^+(u^-)\i$,

(a) are injective and

(b) have the same image.

We denote that common image by $\cb_{>0}$. Note that (b) was proved
in \cite{L94} using the theory of canonical bases. (This is an instance
of $(5)\imp(6)$, see 0.3.)

In \cite{L97} a more direct proof of (b), not relying on the theory
of canonical bases, was given.

From (a),(b) we see that $\Omega$ restricts to an involution
$\cb_{>0}@>>>\cb_{>0}$. This involution is of the form
$c^-\phi'(c^-)\i$ where $\phi':\fU_{>0}@>>>\fU_{>0}$ is an involution.
From the definition  for $u\in\fU_{>0}$ we have
$$\phi'(u)^-B^+(\phi'(u)^-)\i=u^+B^-(u^+)\i.$$

\subhead 7.3\endsubhead
Let $z\in\fU_{>0}$ and let $u\in\fU_{>0}$ so
that $B=u^+B^-(u^+)\i\in\cb_{>0}$. We have
$z^+B(z^+)\i=(zu)^+B^-((zu)^+)\i$,
$$\align&z^-B(z^-)\i=z^-\phi'(u)^-B^+(\phi'(u)^-)\i(z^-)\i\\&
=x^+B^-(x^+)\i=\phi'(x)^-B^+(\phi'(x)^-)\i\endalign$$
where $x\in\fU_{>0}$ satisfies
$\phi'(x)=z\phi'(u)$ that is $x=\phi'(z\phi'(u))$.

Thus if we identify $\fU_{>0}$ and $\cb_{>0}$ by 
$u\lra u^+B^-(u^+)\i$, then the conjugation action of $z^+$
(resp. $z^-$) on $\cb_{>0}$ (with $z\in\fU_{>0}$) becomes the
action $z:u\m zu$  (resp. $z:u\m\phi'(z\phi'(u))$) of $\fU_{>0}$
on $\fU_{>0}$.

\subhead 7.4\endsubhead
We define an involution $\phi:\fU_{>0}@>>>\fU_{>0}$ by
$$\phi=*\phi'*.$$
For $u\in\fU_{>0}$ we have
$$(\phi(u)^-)\i B^+\phi(u)^-=(u^+)\i B^-u^+.$$

\subhead 7.5\endsubhead
Now $G_{>0}$ is equal to the set of all $g\in G(\CC)$ such that
for any irreducible (rational) representation
$\rho$ of $G(\CC)$, the entries of the matrix of
$g$ with respect to the canonical basis of $\rho$
are all in $\RR_{>0}$. (A closely related statement is proved in
\cite{L97, 5.5}.) This is an instance of $(5)\imp(6)$, see 0.3.
A (conjectural) statement in the opposite direction appears in
\cite{L23b, A4}.  This is an instance of $(6)\imp(5)$, see 0.3.

One can show that in type $A$ (so that
$G(\CC)$ is the adjoint group of $GL_n(\CC)$), $G_{>0}$ is the image
under
$$GL_n(\CC)@>>>G(\CC)$$
of the semigroup of strictly totally positive matrices
(matrices with all minors in $\RR_{>0}$) 
studied by Schoenberg, Gantmacher-Krein, Whitney, Loewner.

\subhead 7.6\endsubhead
In \cite{L94} it is shown, using the theory of canonical
bases, that any element of $G_{>0}$ is regular
semisimple when regarded as an element of $G(\CC)$. 
This is an instance of $(5)\imp(6)$, see 0.3.
(In type A this is a result of Gantmacher-Krein; in that case
the theory of canonical bases is not needed.)

In \cite{L94} it is also shown that any $g\in G_{>0}$ is contained
in a unique $B\in\cb_{>0}$.  Hence we have a canonical map $g\m B$,
$\sigma:G_{>0}@>>>\cb_{>0}$. It is known that $\sigma$
is surjective. I
conjecture that each fibre of $\sigma$
is homeomorphic to $\RR_{>0}^\nu$.

\head 8. The semigroup $G(K)$ for a semifield $K$\endhead
\subhead 8.1\endsubhead
Let $K$ be a semifield. We define $G(K)$ to be the semigroup of
permutations of $\fU(K)$ generated by 

(a) the permutations $u\m zu$ for various $z\in\fU(K)$,

(b) the permutations $u\m\phi'_K(z\phi'_K(u))$
for various $z\in\fU(K)$ and

(c) the permutations $x\m S_p(x)$ for various $p:I@>>>K$.

When $K=\RR_{>0}$ we can identify this $G(K)$ with
$G_{>0}$ of \S7.

The permutations (a) define a semigroup imbedding $\fU(K)@>>>G(K)$
with image denoted by $U^+(K)$.
The permutations (b) define a semigroup imbedding $\fU(K)@>>>G(K)$
with image denoted by $U^-(K)$.
The permutations (c) form a commutative subgroup $T(K)$ of $G(K)$.


\widestnumber\key{ABC}
\Refs
\ref\key{A}\by J.Antor\paper Canonical bases via pairing monomials
\jour arxiv:2308.16254\endref
\ref\key{Ca}\by \'E.Cartan\book Oeuvres completes\publ Springer
Verlag\yr1984\endref
\ref\key{Co}\by A.J.Coleman\paper The greatest mathematical paper
of all time\jour Mathematical Intelligencer\vol11\yr1989\pages
29-30\endref
\ref\key{K}\by W.Killing\paper Die Zusammensetzung der stetigen
endlichen Transformationsgruppen\jour Mathematische Ann.
\vol31\yr 1888-1890\endref
\ref\key{G}\by P.Gabriel\paper Unzerlegbare darstellungen\jour Manuscripta Math.\vol6\yr1972\pages71-103\endref
\ref\key{L89}\by G.Lusztig\paper
Modular representations and quantum groups\jour
Contemp.Math.\vol82\yr1989\pages 59-77\endref
\ref\key{L90a}\by G.Lusztig\paper
Finite dimensional Hopf algebras arising from quantized universal
enveloping algebras\jour J. Amer. Math. Soc.\vol3\yr1990\pages
257-296\endref
\ref\key{L90}\by G.Lusztig\paper
 Canonical bases arising from quantized enveloping algebras\jour
J. Amer. Math. Soc.\vol3 \yr1990\pages 447-498\endref
\ref\key{L90b}\by G.Lusztig\paper
 Canonical bases arising from quantized enveloping algebras, II
\inbook Common trends in mathematics and quantum field theories
\bookinfo Progr.of Theor. Phys. Suppl. 102, ed.T.Eguchi et al.\yr1990
\pages175-201\endref
\ref\key{L93}\by G.Lusztig\book Introduction to quantum groups \bookinfo Progr.in Math. 110\publ
Birkh\"auser\publaddr Boston\yr1993\endref
\ref\key{L94}\by G.Lusztig\paper
Total positivity in reductive groups\inbook
Lie theory and geometry\bookinfo Progr. in Math. 123
\publ Birkh\"auser Boston \yr1994\pages531-568\endref
\ref\key{L97}\by G.Lusztig\paper Total positivity and canonical
bases\inbook Algebraic groups and Lie groups, ed. G. I. Lehrer
\publ Cambridge U.Press \yr1997\pages281-295\endref
\ref\key{L23a}\by G.Lusztig\paper The quantum group
$\dot U$ and flag manifolds over the semifield $Z$
\jour Bull. Inst. Math. Acad. Sin.\vol 18\yr2023\pages
235-267\endref
\ref\key{L23b}\by G.Lusztig\paper From Weyl groups to semisimple groups\jour Repres.Th.\vol27\yr2023\pages
51-61\endref
\ref\key{L24}\by G.Lusztig\paper Half circles on flag manifolds over
a semifield, Bull. Inst. Math. Acad. Sin.\vol 19\yr2024\pages1-13
\endref 
\ref\key{L25}\by G.Lusztig\paper History of the canonical basis
and the crystal basis\jour arxiv:2507.20816\endref
\ref\key{R}\by C.Ringel\paper Hall algebras and quantum groups\jour
Invent.Math.\yr1990\vol101\pages583-592\endref
\ref\key{X}\by N.Xi\paper Canonical basis for type $A_3$\jour
Comm. in Alg.\vol27\yr1999\pages5703-5710\endref
\endRefs
\enddocument